\documentclass[dvips,11pt]{article}

\usepackage{graphicx}      
\usepackage{epsfig}        
\usepackage{amsmath,amsthm,amsfonts,latexsym}

\newcommand{\Spec}{\mathrm{Spec}}
\newcommand{\Ass}{\mathrm{Ass}}

\newtheorem{theorem}{Theorem}[section]
\newtheorem{lemma}[theorem]{Lemma}
\newtheorem{proposition}[theorem]{Proposition}
\newtheorem{corollary}[theorem]{Corollary}
\newtheorem{definition}[theorem]{Definition}

\newcounter{counter}

\begin{document}

\title{\normalsize{\bf{CHAINS OF UNUSUAL EXCELLENT LOCAL RINGS}}}

\author{\it{\small{Kai Chen}}}
\date{}

\maketitle

\begin{abstract}
{\noindent Let $(T,M)$ be a complete local domain containing the
integers.
Let $p_{1} \subseteq p_{2} \subseteq \cdots \subseteq p_{n}$ be a chain
of nonmaximal prime ideals of $T$ such that $T_{p_{n}}$ is a regular local
ring. We construct a chain of excellent local
domains $A_{n} \subseteq A_{n-1} \subseteq \cdots \subseteq A_{1}$ such that
for each $1 \leq i \leq n$, the completion of $A_{i}$ is $T$, the 
generic formal fiber of $A_{i}$ is local with maximal ideal $p_{i}$,
and if $I$ is a nonzero ideal of $A_{i}$ then $A_{i}/I$ is complete.
Consequently, if in addition $T$ is a UFD, we can construct a chain of
excellent local UFDs $A_{n} \subseteq A_{n-1} \subseteq \cdots \subseteq
A_{1}$ satisfying the same conditions.
}
\end{abstract}

\section{\textbf{\sc{Introduction}}}

\footnotetext[1]
{2000 {\it Mathematics Subject Classification.} Primary 13J05, 13J10.}
\footnotetext[2]
{{\it Key words and phrases.} Local rings, excellent rings,
completions, generic formal fiber ring.
}

We begin with some basic definitions.
Let $A$ be a local ring with maximal ideal $M$. We use $\hat{A}$ to 
denote the completion of $A$ in the $M$-adic topology. 
If $P$ is a prime ideal of $A$, then the formal fiber ring of $A$ at
$P$ is defined to be $\hat{A} \otimes_A k(P)$, where $k(P) =
A_P/PA_P$. When $A$ is 
a domain, we refer to the formal fiber ring at $(0)$ as the generic
formal fiber ring of $A$. If $p \in \Spec
\hat{A}$ and $K$ is the quotient field of $A$, then when $\hat{A}
\otimes_A K$ is local with maximal ideal $p \otimes_A K$, we say
that $A$ 
has a local generic formal fiber with maximal ideal $p$. 
It is worthwhile to note that it is unusual for an integral domain $A$
to have a local
generic formal fiber, and normally there exist 
many nonzero ideals $I$ of $A$ such that $A/I$ is not complete.
In this paper, we construct chains of excellent rings that do not
satisfy these ``usual'' conditions.

In \cite{LR02}, Loepp and Rotthaus showed that for $T$ a complete local
domain containing the integers with maximal ideal $M$ such that $T/M$
is at least the cardinality of the reals, and $p$ a nonmaximal prime
ideal of $T$ such that $T_p$ is a regular local ring,
there exists an
excellent local domain $A$ such that the completion of $A$ is $T$, the
generic formal fiber of $A$ is local with maximal ideal $p$ and for
any nonzero ideal $I$ of $A$, $A/I$ is complete.

In this paper, we improve the result from \cite{LR02} in two ways. 
We are first able to eliminate the condition that $T/M$ has to have
the cardinality at least that of the reals. In addition, we extend
the result to a chain of excellent local rings. In particular, we show
that if $T$ is a complete local domain containing the integers and $p_1
\subseteq p_2 \subseteq \cdots \subseteq p_n$ a chain of nonmaximal
prime ideals of $T$ such that $T_{p_n}$ is a regular local ring, then
there exists a chain of excellent local domains $A_n \subseteq A_{n-1}
\subseteq \cdots \subseteq A_1$ such that for each $i$, the completion
of $A_i$ is $T$, the generic formal fiber of $A_i$ is local with
maximal ideal $p_i$, and if $I$ is a nonzero ideal of $A_i$ then
$A_i/I$ is complete. As a corollary, we show that if, in addition, $T$
is a UFD, then we may construct a chain of UFDs that satisfy the same
conditions.

The construction of each $A_i$ is similar to that in \cite{LR02}. In
order to guarantee that the completion of each $A_i$ is $T$, we make
use of the result from \cite{H94} (Proposition \ref{p1} below), and
show that $IT \cap A_i=I$ for every finitely generated ideal $I$ of
$A_i$ and $A_i \rightarrow T/M^2$ is onto.  
We also show that $p_i \cap A_i = (0)$ and if $q_i$ is a prime ideal of
$T$ not contained in $p_i$, then $q_i \cap A_i \not= (0)$. This
gives that the generic formal fiber of $A_i$ is local with maximal
ideal $p_i$. We also need to ensure that $A_i/I$ is complete for each
nonzero ideal $I$ of $A_i$.

To extend the result to a chain of excellent local rings, we take an
approach similar to that in \cite{SMALL98}, and break the proof
into two parts. First we construct a chain such that the conditions
hold for the first ring of the chain. Then we use induction to refine
the rest of the chain to obtain the desired rings.

In the construction, we maintain a chain of subrings of $T$
that satisfy some ``nice'' properties. These are similar to the
$p$-subrings in \cite{LR02}. But in order to get rid of the condition
that $T/M$ has to have at least the cardinality of the reals, we take
care to ensure that each subring in the chain is a {\it Small $C$
Avoiding} ($SCA$) subring as defined in \cite{CL03}. We strengthen a
Lemma from \cite{CL03} and use $SCA$-subrings in place of $p$-subrings
in all of our proofs.

Throughout this paper, all rings are commutative with unity. When we
say a ring is local, Noetherian is implied. And we call a ring with
one maximal ideal
that is not necessarily Noetherian quasi-local. We use $(T,M)$ to
denote a quasi-local ring $T$ with maximal ideal $M$, and we use $c$
to denote the cardinality of the reals.

\section{\textbf{\sc{The Construction}}}

We use the following proposition from \cite{H94} to ensure the the
completion of each $A_i$ is $T$.

\begin{proposition}
\label{p1}
If $(R, M \cap R)$ is a quasi-local subring of a complete local ring
$(T, M)$, the map $R \longrightarrow T/M^{2}$ is onto and $IT \cap R = I$ 
for every finitely generated ideal $I$ of $R$, then R is Noetherian and the
natural homomorphism $\hat{R} \longrightarrow T$ is an isomorphism.
\end{proposition}

In our construction, we build subrings of $T$ satisfying some ``nice''
properties. We follow the definition in \cite{CL03} and call them
$SCA$-subrings. 

\begin{definition}
Let $(T, M)$ be a complete local ring and $C$ a set of prime ideals of
$T$. Suppose that $(R, R \cap R)$ is a quasi-local
subring of $T$ such that $|R|<|T|$ and $R \cap P = (0)$ for every $P
\in C$. Then we call $R$ a small $C$-avoiding subring of $T$ and will
denote it by $SCA$-subring.
\end{definition}

Lemma \ref{l7} is well-known, and a proof can be found in \cite{CL03}.
We implicitly use this Lemma in a couple of cardinality arguments.

\begin{lemma}
\label{l7}
Let $(T,M)$ be a complete local ring of dimension at least one. Let
$P$ be a nonmaximal prime ideal of $T$. Then $|T/P| = |T| \geq c$.
\end{lemma}

The following lemma from \cite{CL03} will help us find transcendental
elements. 

\begin{lemma}
\label{l1}
Let $(T, M)$ be a complete local ring such that $\dim T \geq 1$, $C$ a
finite set of nonmaximal prime ideals such that no ideal in $C$ is
contained in another ideal of $C$, and $D$ a subset of $T$ such that
$|D|<|T|$. Let $I$ be an ideal of $T$ such that $I \not\subseteq P$ for
all $P \in C$. Then $I \not\subseteq \bigcup
\lbrace r + P | r \in D, P \in C \rbrace$.
\end{lemma}

We now prove a stronger version of the above lemma which will be
immediately applicable to
our construction. Here we weaken the conditions on $C$ and only
require it to have finitely many maximal elements.

\begin{corollary}
\label{c9}
Let $(T,M)$ be a complete local ring such that $\dim T \geq 1$, $C$ a
set of nonmaximal prime ideals with finitely many maximal elements,
and $D$ a subset of $T$ such that $|D|<|T|$. Let $I$ be an ideal of
$T$ such that $I \not\subseteq P$ for all $P \in C$. Then $I
\not\subseteq \bigcup \{r+P|r \in D, P \in C\}$.
\end{corollary}

\begin{proof}
Let $C'$ be the set of all maximal elements of $C$. So $C'$ is finite,
and no ideal in $C'$ is contained in another ideal of $C'$. Let $I$ be
an ideal of $T$ such that $I \not\subseteq P$ for all $P \in
C$. So $I \not\subseteq P'$ for all $P' \in C'$. Then by Lemma
\ref{l1}, $I \not\subseteq \bigcup \{r+P'| r \in D, P' \in C'\}$.
Now assume, for a contradiction, that $I \subseteq \bigcup \{r+P | r\in
D, P \in C\}$. Then if $x \in I$, we have $x \in r+P$, for some $r \in D,
P \in C$. But $P \subseteq P'$ for some $P' \in C'$, so $x \in r+P'$.
Therefore, $I \subseteq \bigcup \{r+P'| r \in D, P' \in C'\}$, a
contradiction. It follows that $I
\not\subseteq \bigcup \{r+P|r \in D, P \in C\}$.
\end{proof}

Recall that to satisfy the hypothesis of Proposition \ref{p1}, we need
the map $A_i \rightarrow T/M^2$ to be onto. We also want that $A_i/I$
is complete for every nonzero ideal $I$ of $A_i$. Note that if $T$
is the completion of $A_i$, $I$ is a nonzero ideal of $A_i$ and
$A_i \rightarrow T/IT$ is onto, then $\ker (A_i \rightarrow T/IT)= A_i
\cap IT = I$. So 
$\frac{A_i}{I} \cong \frac{T}{IT} \cong \frac{\hat{A_i}}{I\hat{A_i}}
\cong \widehat{\frac{A_i}{I}}$ is complete. To satisfy
this condition and ensure that $A_i \rightarrow T/M^2$ is onto, we
construct $A_i$ such that $A_i \rightarrow T/J_i$ is onto for every ideal
$J_i$ of $T$ that is not contained in the ideal $p_i$. It turns out that
this condition will also help us to show that $A_i$ is excellent.

For $J_i$ an ideal of $T$ not contained in $p_i$, and $\overline{u}
\in T/J_i$, we use the following 2 lemmas  to construct
$SC_iA$-subrings of $T$ satisfying certain properties. These lemmas will
be used to guarantee that $A_i \rightarrow T/J_i$ is onto.

Lemma \ref{l2} is our first attempt to generalize Lemma 3 in
\cite{LR02} to the chain version. Here we show that the desired
properties hold for ``one end'' of the chain. Later we further extend
this result by induction. This approach is similar to that
in \cite{SMALL98}.

\begin{lemma}
\label{l2}
Let $(T, M)$ be a complete local ring with $\dim T \geq 1$, 
$p_{1} \subseteq p_{2} \subseteq \cdots \subseteq p_{n}$ a chain of nonmaximal
prime ideals of $T$ and $J$ an ideal of $T$ with $J \not\subseteq p_{n}$, 
and $J \not\subseteq Q, \forall Q \in \Ass T$. Let $C_i = \Ass T \cup
\{p_i\}$. Let  
$R_{n} \subseteq R_{n-1} \subseteq \cdots \subseteq R_{1}$ be a chain of subrings
of $T$, where $R_{i}$ is an $SC_iA$-subring for each $1 \leq i \leq n$. Let
$\overline{u} \in T/J$. Then there exists a chain of infinite subrings 
$S_{n} \subseteq S_{n-1} \subseteq \cdots \subseteq S_{1}$ of $T$ with
$R_{i} \subseteq S_{i} \subseteq T$, where $S_{i}$ is an $SC_iA$-subring
for each $1 \leq i \leq n$, and
$\overline{u} \in Img(S_{n} \longrightarrow T/J)$. Moreover, if $\overline{u} = 
\overline{0}$, then $S_{n} \cap J \not= (0)$.
\end{lemma}

\begin{proof}
Let $C=\Ass T \cup \{p_{1}, \ldots, p_{n}\}$ and let $P \in C$.
Define $D_{(P)}$ to be a full set of coset representatives of the cosets
$t + P$ that make $(u+t)+P$ algebraic over $\frac{R_{1}}{R_1 \cap P}$
as an element 
of $T/P$. Then $|D_{(P)}| \leq \max(|R_{1}|,\aleph_{0})$. Since
$R_{1}$ is an $SC_1A$-subring, we have $|R_{1}| < |T|$. It follows
that $|D_{(P)}| < |T|$. Clearly, $J \not\subseteq P, \forall P \in
C$. 
Let $D = \bigcup_{P \in C}D_{(P)}$. Note that $|D| < |T|$.
Use Corollary {\ref{c9}} to find $x \in J$ such that 
$x \not\in \bigcup\{r+P|r \in D, P \in C\}$. We define
$S_{i} = R_{i}[x+u]$ localized at $R_{i}[x+u] \cap M$.
Since we are adjoining the same element, clearly 
$S_{n} \subseteq S_{n-1} \subseteq \cdots \subseteq S_{1}$. 
And note that each $S_i$ is infinite.
We claim the $S_{i}$'s are the desired $SC_iA$-subrings.

Since $|R_{i}| < |T|$ and $S_{i} = max(|R_{i}|, \aleph_{0})$,
we have $|S_{i}| < |T|$. Suppose $P_i \in C_i$. 
Then if $f_{i} \in R_{i}[u+x] \cap P_i$, we have 
$$
f_{i} = r_{i,n}(u+x)^n + \cdots + r_{i,1}(u+x) + r_{i,0} \in P_i.
$$
So $f_{i} \equiv 0$ modulo $P_i$. Since 
$(u+x)+P_i$ is transcendental over $\frac{R_{1}}{R_1 \cap P_i}$, it is
transcendental over $\frac{R_{i}}{R_i \cap P_i} \cong R_i$. Thus
$r_{i,j} \equiv 0$ modulo $P_i$ for each $j$. Then 
$r_{i,j} \in P_i$, and $r_{i,j} \in P_i \cap R_{i} = (0)$. 
Hence, $S_{i} \cap P_i = (0)$ and it follows that 
$S_{i}$ is a $SC_iA$-subring.
Note that under the map $S_{n} \longrightarrow T/J$,
$u+x$ is mapped to $u+J$. Hence 
$\overline{u} \in Img(S_{n} \longrightarrow T/J)$.
Also, if $\overline{u} = \overline{0}$, then $u+x \in J$. 
Since $(u+x)+p_{n}$ is transcendental over $R_{n}$ as an element
of $T/p_{n}$, we have $u+x \not= 0$. It follows that 
$S_{n} \cap J \not= (0)$.
\end{proof}

Now we are ready to give the full chain version of Lemma 3 in
\cite{LR02}. We argue inductively from Lemma \ref{l2}.

\begin{lemma}
\label{l3}
Let $(T, M)$ be a complete local ring with $\dim T \geq 1$, 
$p_{1} \subseteq p_{2} \subseteq \cdots \subseteq p_{n}$ a chain of
nonmaximal prime ideals of $T$ and $J_{1}$, $J_{2}$ $\ldots$ $J_{n}$  
$n$ ideals of $T$ with $J_{i} \not\subseteq p_{i}$ for each $1 \leq i \leq n$, 
and $J_i \not\subseteq Q, \forall Q \in \Ass T$. Let $C_i = \Ass T
\cup \{p_i\}$. Let $R_{n} \subseteq R_{n-1} \subseteq \cdots \subseteq
R_{1}$ be a chain of subrings 
of $T$, where $R_{i}$ is an $SC_iA$-subring for each $1 \leq i \leq n$. Let
$\overline{u}_{i} \in T/J_{i}$ for each $i$. Then there exists a chain
of infinite subrings  
$S_{n} \subseteq S_{n-1} \subseteq \cdots \subseteq S_{1}$ of $T$ with
$R_{i} \subseteq S_{i} \subseteq T$, where $S_{i}$ an $SC_iA$-subring 
for each $1 \leq i \leq n$, and $\overline{u}_{i} \in Img(S_{i}
\longrightarrow T/J_{i})$ for every $i$. Moreover, if 
$\overline{u}_{i} = \overline{0}$, then $S_{i} \cap J_{i} \not= (0)$.
\end{lemma}

\begin{proof}
We will induct on $n$. The base case is clear if we let $n=1$ in Lemma
{\ref{l2}}. Assume given a chain of $n-1$ prime ideals we can find a
chain of  
$n-1$ desired $SC_iA$-subrings. Now consider the case with $n$ prime
ideals. Use Lemma {\ref{l2}} to find a chain of $n$ subrings
$S'_{n} \subseteq S'_{n-1} \subseteq \cdots \subseteq S'_{1}$ such that
$S'_{i}$ is an $SC_iA$-subring, $R_{i} \subseteq S'_{i} \subseteq T$, 
$\overline{u}_{n} \in Img(S_{n} \longrightarrow T/J_{n})$, and
if $\overline{u}_{n} = \overline{0}$, then 
$S_{n} \cap J_{n} \not= (0)$.

Next use the inductive hypothesis on the chain 
$S'_{n-1} \subseteq S'_{n-2} \subseteq \cdots \subseteq S'_{1}$
to find
$S_{n-1} \subseteq S_{n-2} \subseteq \cdots \subseteq S_{1}$ such that
for $1 \leq i \leq n-1$, 
$S_{i}$ is an $SC_iA$-subring, 
$R_{i} \subseteq S'_{i} \subseteq S_{i} \subseteq T$, 
$\overline{u}_{i} \in Img(S_{i} \longrightarrow T/J_{i})$, and
if $\overline{u}_{i} = \overline{0}$, then 
$S_{i} \cap J_{i} \not= (0)$.

Finally, let $S_{n} = S'_{n}$. We have
$S_{n} = S'_{n} \subseteq S'_{n-1} \subseteq S_{n-1}$.
Hence, we have constructed the chain:
$S_{n} \subseteq S_{n-1} \subseteq \cdots \subseteq S_{1}$ as desired.
By induction the lemma holds.
\end{proof}

The following two lemmas are the
generalization of Lemma 2.6 from \cite{CL03}. We need this to show that
$I_iT \cap A_i = I_i$ for each finitely generated ideal $I_i$ of
$A_i$. The proof again is broken up into two steps similar to the ones in
the previous two lemmas.

\begin{lemma}
\label{l8}
Let $(T, M)$ be a complete local ring with $M \not\in \Ass T$
and $\dim T \geq 1$. Let 
$p_{1} \subseteq p_{2} \subseteq \cdots \subseteq p_{n}$ be a chain of nonmaximal
prime ideals of $T$. Let $C_i = \Ass T \cup \{p_i\}$
and $R_{n} \subseteq R_{n-1} \subseteq \cdots \subseteq R_{1}$
a chain of subrings of $T$ where $R_{i}$ is an $SC_iA$-subring 
for each $1 \leq i \leq n$.
Suppose $I$ is a finitely generated ideal of $R_n$ and $c \in IT \cap R_n$. 
Then there exists a chain of subrings
$S_{n} \subseteq S_{n-1} \subseteq \cdots \subseteq S_{1}$ of $T$ such that for each
$1 \leq i \leq n$, $R_{i} \subseteq S_{i} \subseteq T$, $S_{i}$ is an $SC_iA$-subring
and $c \in IS_{n}$.
\end{lemma}

\begin{proof}
We will induct on the number of generators of $I$. If $I$ is
principal, $I = aR_n$, for some $a \in R_n$. If $a=0$, then $I=(0)$
and $c \in IT \cap R_n$ implies $c=0$. So $S_i = R_i$ gives the desired
$SC_iA$-subrings. So we consider the case when $a \not= 0$. In this
case, $c = au$ for some $u \in T$. We claim that $S_i =
R_i[u]_{(R_i[u]\cap M)}$ are the desired $SC_iA$-subrings. To see
this, first note that $|S_i|<|T|$ for each $i = 1,2,\ldots,n$. Now let
$P_i \in C_i$ and suppose $f_i = r_nu^n + \cdots + r_1u+r_0 \in R_i[u]
\cap P_i$. Multiplying through by $a^n$, we obtain
$a^nf_i = r_n(au)^n + \cdots + r_1a^{n-1}(au)+r_0a^n$. 
Since $a \in R_n \subseteq R_i$, and $R_i \cap P_i = (0)$ for each
$P_i \in C_i$, it follows that
$a^nf_i = r_nc^n + \cdots + r_1a^{n-1}c+r_0a^n \in P_i \cap R_i =
(0).$ Since $C_n$ contains all associated primes of $T$, and $a \in
R_n$ which is a $SC_nA$-subring, we know $a$
is not a zero divisor in $T$. Then it must be the case that $f=0$. 
Thus $S_i$ is an $SC_iA$-subring and this proves the base case.

Now suppose the Lemma holds for all $I$ generated by $m-1$
elements. Consider the case when $I$ has $m$ generators. 
Let $I = (y_1,\ldots,y_m)R_n$. Since $c \in IT \cap R_n$,
$c = y_1t_1 + \cdots + y_mt_m$ for some $t_1, \ldots, t_m \in T$.
Note that by adding $0$, for any $t \in T$,  
$c = y_1t_1 + y_1y_2t-y_1y_2t+y_2t_2
+ \cdots + y_mt_m = y_1(t_1+y_2t) + y_2(t_2-y_1t)+y_3t_3 + \cdots
y_mt_m$. Let $x_1 = t_1 + y_2t$ and $x_2 = t_2 - y_1t$, where we will
choose $t$ later. Let $C =\Ass T \cup \{p_1,\ldots,p_n\}$ and $P \in 
C$. Note that $y_2 \not\in P$ since $y_2 \in R_n$, $y_2 \not=0$, and
$R_n \cap P = (0)$. Thus, $t+P\not=t'+P$ implies $(t_1+y_2t)+P \not=
(t_1+y_2t')+P$. Let $D_{(P)}$ be a full set of coset representatives
of the cosets $t+P$ that makes $x_1+P$ algebraic over $\frac{R_1}{R_1
  \cap P}$. Let $D =
\bigcup_{P \in C} D_{(P)}$. Note that $|D| < |T|$. We use Corollary
\ref{c9} with $I=T$ to find an element $t\in T$ such that $x_1+P$ is
transcendental over $\frac{R_1}{R_1 \cap P}$ for every $P \in C$. Then
$x_1+P_i$ is transcendental over $\frac{R_i}{R_i \cap P_i} \cong R_i$
for every $P_i \in C_i$. 
As in the proof of Lemma \ref{l2}, we may show
that $R'_i = R_i[x_1]_{(R_i[x_1]\cap M)}$ are $SC_iA$-subrings.
Let $I' = (y_2,\ldots,y_m)R'_n$ and $c^* = c-y_1x_1$. So $c^* \in I'T
\cap R'_n$. Then by our inductive hypothesis, there exists a chain of
subrings $S_n \subseteq \cdots \subseteq S_1$ such that $R'_i
\subseteq S_i \subseteq T$, $S_i$ is an $SC_iA$-subring, and $c^* \in
I'S_n$. Thus $c^* = y_2s_2 + \cdots + y_ms_m$ for $s_1,\ldots,s_m \in
S_n$. It follows that $c=y_1x_1+y_2s_2+\cdots+y_ms_m \in IS_n$. And since
$R_i \subseteq R'_i$, we conclude that $R_i \subseteq S_i \subseteq T$
and the $S_i$'s are the desired $SC_iA$-subrings.
\end{proof}

\begin{lemma}
\label{l4}
Let $(T, M)$ be a complete local ring with $M \not\in \Ass T$
and $\dim T \geq 1$. Let 
$p_{1} \subseteq p_{2} \subseteq \cdots \subseteq p_{n}$ be a chain of
nonmaximal prime ideals of $T$. Let $C_i = \Ass T \cup \{p_i\}$
and $R_{n} \subseteq R_{n-1} \subseteq \cdots \subseteq R_{1}$
a chain of subrings of $T$ where $R_{i}$ is an $SC_iA$-subring 
for each $1 \leq i \leq n$.
Let $I_{1}$, $I_{2}$ $\ldots$ $I_{n}$ and $c_{1}$, $c_{2}$ $\ldots$ $c_{n}$
be such that $I_{i}$ is a finitely generated ideal of $R_{i}$ with
$c_{i} \in I_{i}T \cap R_{i}$. Then there exists a chain of subrings
$S_{n} \subseteq S_{n-1} \subseteq \cdots \subseteq S_{1}$ of $T$ such that
for each $1 \leq i \leq n$, $R_{i} \subseteq S_{i} \subseteq T$,
$S_{i}$ is an $SC_iA$-subring and $c_{i} \in I_{i}S_{i}$.
\end{lemma}

\begin{proof}
We will induct on $n$. The base case clearly holds if we let $n=1$ in
Lemma \ref{l8}. Assume that given a chain of $n-1$ prime ideals and
the corresponding chain of subrings, we can find a chain of $n-1$
desired $SC_iA$-subrings. Now consider the case with $n$ prime
ideals. Use Lemma \ref{l8} to construct a chain of subrings
$S'_n \subseteq S'_{n-1} \subseteq \cdots \subseteq S'_1$ such that
for each $i$, $R_i \subseteq S'_i \subseteq T$, $S'_i$ is a
$SC_iA$-subring, and $c_n \in I_nS'_n$. We let $S_n = S'_n$ and use
our inductive hypothesis on the chain $S'_{n-1} \subseteq \cdots
\subseteq S'_1$ to find a chain of subrings $S_{n-1} \subseteq \cdots
\subseteq S_1$ such that $R_i \subseteq S'_i \subseteq S_i \subseteq
T$, $S_i$ is $SC_iA$-subring, and $c_i \in I_iS_i$. Since $S_n = S'_n
\subseteq S'_{n-1} \subseteq S_{n-1}$, we have constructed a chain of
subrings $S_n \subseteq \cdots \subseteq S_1$ satisfying all the
conditions. By induction, our lemma holds.
\end{proof} 

The following definition is from \cite{H93}.

\begin{definition}
Let $\Omega$ be a well-ordered set and $\alpha \in \Omega$. We define
$\gamma (\alpha) = sup \lbrace \beta \in \Omega | \beta < \alpha \rbrace$.
\end{definition}

Lemma \ref{l5} is the generalization of Lemma 5 from \cite{LR02}. Here
we construct a chain of $SC_iA$-subrings that simultaneously satisfy many
of the desired properties. Condition ($iii$) and ($v$) will help us to
show that the completion of $A_i$ is $T$. Condition ($iv$) is needed to
ensure that the generic formal fiber of $A_i$ is local with maximal
ideal $p_i$.

\begin{lemma}
\label{l5}
Let $(T, M)$ be a complete local ring with $\dim T \geq 1$, 
$p_{1} \subseteq p_{2} \subseteq \cdots \subseteq p_{n}$ a chain of nonmaximal
prime ideals of $T$ and $J_{1}$, $J_{2}$ $\ldots$ $J_{n}$  
$n$ ideals of $T$ with $J_{i} \not\subseteq p_{i}$ for each $1 \leq i \leq n$. 
Let $C_i = \Ass T \cup \{p_i\}$. Let 
$R_{n} \subseteq R_{n-1} \subseteq \cdots \subseteq R_{1}$ be a chain of subrings
of $T$, where $R_{i}$ is an $SC_iA$-subring. Let $\overline{u}_{i} \in
T/J_{i}$ for every $i$. 
Then there exists a chain of subrings 
$S_{n} \subseteq S_{n-1} \subseteq \cdots \subseteq S_{1}$ of $T$ such that
for every $1 \leq i \leq n$
\begin{list}{(\roman{counter})}
{\usecounter{counter} \setlength{\rightmargin}{\leftmargin}}
\item $S_{i}$ is an $SC_iA$-subring, 
\item $R_{i} \subseteq S_{i} \subseteq T$, 
\item $\overline{u}_{i} \in Img(S_{i} \longrightarrow T/J_{i})$, 
\item if $\overline{u}_{i} = \overline{0}$, then $S_{i} \cap J_{i} \not= (0)$,
and
\item For every finitely generated ideal $I$ of $S_{i}$, we have
$IT \cap S_{i} = I$.
\end{list}
\end{lemma}

\begin{proof}
Use Lemma {\ref{l3}} to find a chain of infinite subrings
$R_{n,0} \subseteq R_{n-1,0} \subseteq \cdots \subseteq R_{1,0}$ such that
for each $1 \leq i \leq n$, $R_{i,0}$ is an $SC_iA$-subring,
$R_{i} \subseteq R_{i,0} \subseteq T$, $\overline{u}_{i} \in 
Img(R_{i,0} \longrightarrow T/J_{i})$ and if $\overline{u}_{i} = 
\overline{0}$, then $R_{i,0} \cap J_{i} \not= (0)$.
We construct $S_{i}$ to contain $R_{i,0}$, so conditions $(ii)-(iv)$
will follow automatically. Now for each
$1 \leq i \leq n$,  define
$$
\Omega_{i} = \{(I_{i}, c_{i}) | I_{i} \mbox{ is a finitely generated ideal
of } R_{i,0} \mbox{ and } c_{i} \in I_{i}T \cap R_{i,0}\}.
$$
Since $I_{i}$ can be $R_{i,0}$, we have $|R_{i,0}| \leq |\Omega_{i}|$,
and since $|R_{i,0}|$ is infinite and there are $|R_{i,0}|$ finite
subsets of $R_{i,0}$, we have 
$|\Omega_{i}| \leq |R_{i,0}|$. So $|\Omega_{i}| = |R_{i,0}|$.
Well-order each $\Omega_{i}$ so that it does not have a maximal
element.
Note that $((0),0) \in \Omega_i$ for all $i$. 
Suppose $\Omega_I$ is a set with the maximal cardinality, let $\Psi$
be its index set. We use $\Psi$ to index each $\Omega_i$, and when
$|\Omega_i| < |\Omega_I|$, we simply append extra $((0),0)$'s.
Now define
$$
\Omega = \Delta\{(\Omega_{1} \times \Omega_{2} \times \cdots 
\times \Omega_{n})\}.
$$
where $\Delta$ denotes the diagonal.
We can naturally well-order $\Omega$ using $\Psi$.
Next, we recursively define $n$ families of $SC_iA$-subrings, starting with
$R_{n,0} \subseteq R_{n-1,0} \subseteq \cdots \subseteq R_{1,0}$. We take care
to preserve the ascending chain structure at each step.
Let $\alpha \in \Omega$. Assume 
$R_{n,\beta} \subseteq R_{n-1,\beta} \subseteq \cdots \subseteq R_{1,\beta}$
has been defined for each $\beta < \alpha$.
$\gamma(\alpha) = ((I_{1},c_{1}), (I_{2},c_{2}), ..., (I_{n},c_{n}))$ 
for some finitely generated ideal $I_{i}$ of $R_{i,0}$ and 
$c_{i} \in I_{i}T \cap R_{i,0}$.
If $\gamma(\alpha) < \alpha$, use Lemma {\ref{l4}} to find 
$R_{n,\alpha} \subseteq R_{n-1,\alpha} \subseteq \cdots \subseteq
R_{1,\alpha}$ such that $R_{i,\alpha}$ is an $SC_iA$-subring,
$R_{i,\gamma(\alpha)} \subseteq R_{i, \alpha} \subseteq T$
and $c_{i} \in I_{i}R_{i,\alpha}$.
If $\gamma(\alpha) = \alpha$, define 
$R_{i,\alpha} = \bigcup_{\beta<\alpha}R_{i,\beta}$.
We still need to show that the chain condition holds in this case.
Let $1 \leq a < b \leq n$, and $\zeta \in R_{b, \alpha}$. Hence
$\zeta \in R_{b, \beta}$ for some $\beta < \alpha$. Since each of
$R_{i, \beta}$ has been defined to preserve the chain, we have
$\zeta \in R_{a, \beta} \subseteq R_{a, \alpha}$. It follows that
$R_{b,\alpha} \subseteq R_{a,\alpha}$, as desired. Now, define
$R_{i,1} = \bigcup_{\alpha \in \Omega}R_{i,\alpha}$.
Since $R_{n,\alpha} \subseteq R_{n-1,\alpha} \subseteq \cdots \subseteq
R_{1,\alpha}$ holds for each $\alpha \in \Omega$, it follows that
$R_{n,1} \subseteq R_{n-1,1} \subseteq \cdots \subseteq R_{1,1}$.
It is easy to verify that $R_{i,1}$ is an $SC_iA$-subring, for each
$i=1,2,\ldots, n$. 
Also, if $I_{i}$ is a finitely generated ideal of $R_{i,0}$, and
$c_{i} \in I_{i}T \cap R_{i,0}$, then 
$((I_{1}, c_{1}), \ldots, (I_{i}, c_{i}), \ldots (I_{n}, c_{n}))
 = \gamma(\alpha)$
for some $\alpha \in \Omega$ with $\gamma(\alpha) < \alpha$. 
So, $c_{i} \in I_{i}R_{i,\alpha} \subseteq I_{i}R_{i,1}$.
Hence $I_{i}T \cap R_{i,0} \subseteq I_{i}R_{i,1}$.

Repeat this process to obtain a chain of $SC_iA$-subrings
$R_{n,2} \subseteq R_{n-1,2} \subseteq \cdots \subseteq R_{1,2}$ such that
$I_{i}T \cap R_{i,1} \subseteq I_{i}R_{i,2}$
Continue to construct an ascending chain
$R_{i,0} \subseteq R_{i,1} \subseteq \cdots$ such that for each $j$,
$R_{n,j} \subseteq R_{n-1,j} \subseteq \cdots \subseteq R_{1,j}$ is 
a chain of $SC_iA$-subrings,
and for each finitely generated ideal $I_{i}$ of $R_{i,m}$, 
$I_{i}T \cap R_{i,m} \subseteq I_{i}R_{i,m+1}$.
Then $S_{i} = \bigcup_{j=1}^{\infty}R_{i,j}$ is an $SC_iA$-subring.
Clearly $S_{n} \subseteq S_{n-1} \subseteq \cdots \subseteq S_{1}$.
If $I$ is a finitely generated ideal of $S_{i}$, then some $R_{i,n}$
contains a generating set for $I$, say
$\{y_{1}, \ldots, y_{k}\} \subseteq R_{i,n}$. Clearly
$I \subseteq IT \cap S_{i}$. Now let $c \in IT \cap S_{i}$, then 
$c \in R_{i,m}$ for some $m \geq n$. Hence 
$c \in (y_{1}, \ldots, y_{k})R_{m+1} \subseteq I$.
It follows that $I = IT \cap S_{i}$.
\end{proof}

The following lemma is a generalization of Lemma 6 in
\cite{LR02}. Here we construct a chain of rings $A_i$ satisfying all
conditions we desire except that $A_i$ might not be excellent.

\begin{lemma}
\label{l6}
Let $(T, M)$ be a complete local ring with $\dim T \geq 1$ and such that no
integer of $T$ is a zero-divisor. Let 
$p_{1} \subseteq p_{2} \subseteq \cdots \subseteq p_{n}$ be a chain of
nonmaximal prime ideals of $T$, such that $p_{1}$ contains all the
associated prime ideals of $T$. 
Suppose that $p_{n}$ intersected with the prime subring of $T$
is the zero ideal. Then there exists a chain
of local domains $A_{n} \subseteq A_{n-1} \subseteq \cdots \subseteq A_{1}$
such that for each $1 \leq i \leq n$
\begin{list}{(\roman{counter})}
{\usecounter{counter} \setlength{\rightmargin}{\leftmargin}}
\item $\hat{A}_{i} = T$,
\item if $P_{i}$ is a nonzero prime ideal of $A_{i}$, then
$T \otimes_{A_{i}} k(P_{i}) \cong k(P_{i})$,
where $k(P_i)=\frac{A_{P_i}}{P_iA_{P_i}}$,
\item the generic formal fiber of $A_{i}$ is local with maximal ideal
$p_{i}$, and
\item if $I$ is a nonzero ideal of $A_{i}$, then $A_{i}/I$ is
  complete.
\end{list}
\end{lemma}

\begin{proof}
Let $C_i = \Ass T \cup \{p_i\}$. 
Define $n$ sets $\Omega_{1}$, $\Omega_{2}$, $\ldots$, $\Omega_{n}$
to be such that
$$
\Omega_{i} = \lbrace u+J_{i} \in T/J_{i} | J_{i} \mbox{ is an ideal
  of } T \mbox{ with } J_{i} \not\subseteq p_{i}\rbrace.
$$
Since $T$ is Noetherian, then each ideal of $T$ is finitely generated.
Hence $|\{J_{i} | J_{i}$ is ideal of $T$ and $J_{i} \not\subseteq p_{i}\}| 
\leq |T|$. Now, if $J$ is an ideal of $T$, then $|T/J| \leq |T|$.
It follows that for each $1 \leq i \leq n$,
$|\Omega_{i}| \leq |T|$. Well order $\Omega_{i}$ so that each element
has fewer than $|\Omega_{i}|$ predecessors. 
Note that $0+M \in \Omega_i$ for each $i$.
By using the same technique as in the proof of Lemma \ref{l5}, 
we use the same index set for all $n$ 
orderings, call this index set $\Psi$. Next define
$$
\Omega = \Delta\{(\Omega_{1} \times \Omega_{2} \times \cdots 
\times \Omega_{n})\}.
$$
We can naturally well order $\Omega$ using $\Psi$. Let $0$ denote the
first element of $\Omega$. Let $R'_{0}$ be the prime subring of $T$
and for each $1 \leq i \leq n$, let $R_{i,0} = {R'_{0}}_{(R'_{0} \cap
  M)}$.  Note that $R_{i,0}$ is an $SC_iA$-subring and the condition 
$R_{n,0} \subseteq R_{n-1,0} \subseteq \cdots \subseteq R_{1,0}$
holds trivially. 

Next, we recursively define $n$ families of $SC_iA$-subrings and take care
to preserve the chain structure of the $SC_iA$-subrings at each step.
$R_{n,0} \subseteq R_{n-1,0} \subseteq \cdots \subseteq R_{1,0}$
is already defined. 
Let $\lambda \in \Omega$. Assume 
$R_{n,\beta} \subseteq R_{n-1,\beta} \subseteq \cdots \subseteq R_{1,\beta}$
has been defined for each $\beta < \lambda$. Then 
$\gamma(\lambda) = ((u_{1}+J_{1}), (u_{2}+J_{2}), \ldots, (u_{n}+J_{n}))$ 
for some ideals $J_{i}$ of $T$ with $J_{i} \not\subseteq p_{i}$. 
If $\gamma(\lambda) < \lambda$, use Lemma {\ref{l5}} to obtain a chain of
$SC_iA$-subrings 
$R_{n,\lambda} \subseteq R_{n-1,\lambda} \subseteq \cdots \subseteq
R_{1,\lambda}$ 
such that $R_{i,\gamma(\lambda)} \subseteq R_{i,\lambda} \subseteq T$,
and $u_{i}+J_{i} \in Img(R_{i, \lambda} \longrightarrow T/J_{i})$.
Moreover, if $u_{i}+J_{i} = 0+J_{i}$, then $R_{i, \lambda} \cap J_{i}
\not= (0)$, and
for every finitely generated ideal $I$ of $R_{i, \lambda}$, we have
$IT \cap R_{i, \lambda} = I$.
If $\gamma(\lambda) = \lambda$, define 
$R_{i, \lambda} = \bigcup_{\beta < \lambda}R_{i, \beta}$. 
We still need to show the chain condition holds in this case.
Let $1 \leq a < b \leq n$, and $\zeta \in R_{b, \lambda}$. Hence
$\zeta \in R_{b, \beta}$ for some $\beta < \lambda$. Since each of
$R_{i, \beta}$ has been defined to preserve the chain, we have
$\zeta \in R_{a, \beta} \subseteq R_{a, \lambda}$. Thus
$R_{n,\lambda} \subseteq R_{n-1,\lambda} \subseteq \cdots \subseteq R_{1,\lambda}$
holds for each $\lambda \in \Omega$ and $R_{i, \lambda}$ is an $SC_iA$-subring.
Let $A_{i} = \bigcup_{\lambda \in \Omega}R_{i,\lambda}$, clearly
$A_{n} \subseteq A_{n-1} \subseteq \cdots \subseteq A_{1}$ and we claim this is
the chain of the desired domains.

We first prove condition $(iii)$. Since each $R_{i, \lambda}$ is a 
$SC_iA$-subring, we have $R_{i, \lambda} \cap p_{i} = (0)$. 
Hence, $A_{i} \cap p_{i} = (0)$. Next, let $J_{i}$ be an ideal of $T$
with $J_{i} \not\subseteq p_{i}$, then $0+J_{i} \in \Omega_{i}$.
So there exists $\lambda \in \Omega$ with $\gamma(\lambda)<\lambda$
and $\gamma(\lambda) = ((u_{1}+J_{1}), (u_{2}+J_{2}), \ldots, 
(0+J_{i}), \ldots, (u_{n}+J_{n}))$.
By our construction, $R_{i,\lambda} \cap J_{i} \not= (0)$. It follows
that $A_{i} \cap J_{i} \not= (0)$, $\forall J_{i}$ an ideal of $T$ with
$J_{i} \not\subseteq p_{i}$. Thus, the generic formal fiber of $A_{i}$
is local with maximal ideal $p_{i}$.

To show $(i)$, we make use of Proposition {\ref{p1}}. Since $p_{i}$
is nonmaximal, we have $M^{2} \not\subseteq p_{i}$. Thus by our construction,
the map $A_{i} \longrightarrow T/M^{2}$ is onto. 
For every finitely generated ideal $I$ of $A_{i}$, clearly 
$I \subseteq IT \cap A_{i}$. Let $I = (y_{1}, y_{2}, \ldots, y_{k})$
be a finitely generated ideal of $A_i$ and
$c \in IT \cap A_{i}$. We have $(c, y_{1}, \ldots, y_{k}) \subseteq
R_{i, \lambda}$ for some $\lambda \in \Omega$ with 
$\gamma(\lambda) < \lambda$. By construction, 
$(y_{1}, \ldots, y_{k})T \cap R_{i, \lambda} = 
(y_{1}, \ldots, y_{k})R_{i, \lambda}$.
Hence, $c \in (y_{1}, \ldots, y_{k})T \cap R_{i, \lambda}
= (y_{1}, \ldots, y_{k})R_{i, \lambda} \subseteq I$. We thus have
$IT \cap A_{i} = I$. It follows from Proposition {\ref{p1}}
that $A_{i}$ is Noetherian and the completion of $A_{i}$ is $T$.

Next, let $I$ be a nonzero ideal of $A_{i}$. Let $J_{i} = IT$.
Suppose $J_{i} \subseteq p_{i}$, then $I \subseteq J_{i} \cap A_{i}
\subseteq p_{i} \cap A_{i} = (0)$, a contradiction. Hence, 
$J_{i} \not\subseteq p_{i}$. It follows that the map
$A_{i} \longrightarrow T/IT$ is surjective.
The kernel of this map is $A_{i} \cap IT = I$. Hence
$A_{i}/I \cong T/IT$, which implies that $A_{i}/I$ is complete.

Finally, if $P_i$ is a nonzero prime ideal of $A_i$, then by the above
paragraph, we have that
$A_{i}/P_{i} \cong T/P_{i}T$. It then follows that
$T \otimes_{A_{i}} k(P_{i})\cong (T/P_{i}T)_{\overline{A_{i}-P_{i}}}
\cong (A_{i}/P_{i})_{\overline{A_{i}-P_{i}}} \cong
{A_{i}}_{P_{i}}/P_i{A_{i}}_{P_{i}} = k(P_{i})$.

And this completes the proof. 
\end{proof}

Now we prove that each $A_i$ is excellent and thus conclude our main
theorem. 

\begin{theorem}
\label{t1}
Let $(T, M)$ be a complete local domain containing the integers.
Let $p_{1} \subseteq p_{2} \subseteq \cdots \subseteq p_{n}$
be a chain of nonmaximal prime ideals of $T$ such that $T_{p_{n}}$ is
a regular local ring. Then there exists a chain of excellent local
domains $A_{n} \subseteq A_{n-1} \subseteq \cdots  
\subseteq A_{1}$ such that for each $1 \leq i \leq n$, the completion of
$A_{i}$ is $T$, the generic formal fiber of $A_{i}$ is local with maximal
ideal $p_{i}$ and if $I$ is a nonzero ideal of $A_{i}$, then
$A_{i}/I$ is complete.
\end{theorem}

\begin{proof}
Notice that since $T$ is a local domain, in order for the chain $p_1
\subseteq \cdots \subseteq p_n$ to exist, we must have $\dim T \geq 1$.
We use Lemma {\ref{l6}} to construct the chain of local domains
$A_{n} \subseteq A_{n-1} \subseteq \cdots \subseteq A_{1}$.
Then we only need to show that for each $1 \leq i \leq n$, $A_{i}$ is
excellent.
Since $A_{i}$ is local (and thus Noetherian), to show $A_{i}$ is
universally catenary, we only need to show $A_{i}$ is formally 
equidimensional. Hence, it suffices to show that $\hat{A_{i}} = T$
is equidimensional. This is true since $T$ is domain.

Next, we show $A_{i}$ is $G-ring$ for each $i$.
Let $P_{i}$ be a nonzero prime ideal of $A_{i}$. 
By Lemma {\ref{l6}}, $T \otimes_{A_{i}} k(P_{i}) \cong k(P_{i})$
is regular.
Let $L$ be a finite field extension of $k(P_{i})$. We want to show
$T \otimes_{A_{i}} L$ is regular. Here we make use of a fact that
if $M$ is an $R$-module, then $M \otimes_{R} R \cong M$.
Hence we have:
$T \otimes_{A_{i}} L \cong T \otimes_{A_{i}} (k(P_{i})
\otimes_{k(P_{i})} L) = (T \otimes_{A_{i}} k(P_{i}))
\otimes_{k(P_{i})} L \cong k(P_{i}) \otimes_{k(P_{i})} L
\cong L$ is regular.

Now the only thing left to show is that the generic formal fiber
of $A_{i}$ is geometrically regular. 
First, recall that if $R$ is a regular local ring, then $R_P$ is a
regular local ring $\forall P \in \Spec R$. So $T_{p_n}$ being a RLR
implies that $T_{p_i} \cong (T_{p_n})_{p_iT_{p_n}}$ is a RLR. Then,
$T \otimes_{A_{i}} k(0) \cong T_{p_{i}}$ is regular.
Next, we want to show that $T \otimes_{A_{i}} L$ is regular,
$\forall$ finite field extension $L$ of $k(0)$.
Since $k(0) = {A_{i}}_{(0)}/(0){A_{i}}_{(0)}$ is the quotient field of 
$A_{i}$, it is a field of characteristic zero. Thus by Rotthaus'
notes in \cite{Ro97}, we may assume that $L$ is purely inseparable.
Since $k(0)$ is characteristic zero, it must be that $L = k(0)$.
Hence, it follows that $A_{i}$ is excellent for each 
$1 \leq i \leq n$.
\end{proof}

Finally, the following corollary extends the results to UFDs.

\begin{corollary}
\label{c3}
Let $(T, M)$ be a complete local unique factorization domain 
containing the integers.
Let $p_{1} \subseteq p_{2} \subseteq \cdots \subseteq p_{n}$
be a chain of nonmaximal prime ideals of $T$ such that $T_{p_{n}}$ is
a regular local ring. Then there exists
a chain of excellent local UFDs $A_{n} \subseteq A_{n-1} \subseteq \cdots 
\subseteq A_{1}$ such that for each $1 \leq i \leq n$, the completion of
$A_{i}$ is $T$, the generic formal fiber of $A_{i}$ is local with maximal
ideal $p_{i}$ and if $I$ is a nonzero ideal of $A_{i}$, then
$A_{i}/I$ is complete.
\end{corollary}

\begin{proof}
We use Theorem \ref{t1} to construct a chain of excellent local
domains $A_n \subseteq A_{n-1} \subseteq \cdots \subseteq A_1$. It remains
to show that each $A_i$ is a UFD. But this follows immediately from
the fact (see Exercise 20.4 in \cite{Ma86}) that for a local ring $A$,
if $\hat{A}$ is a UFD, then so is $A$. 
\end{proof}

\vspace{0.6cm}
\begin{center}
{\bf Acknowledgements}
\end{center}
I thank S. Loepp for her helpful advice and
editorial comments. I also thank D. Jensen for pointing out a
useful fact which immediately leads to Corollary \ref{c3}.

\end{document}